\documentclass[11pt]{amsart}
\usepackage{amssymb,mathrsfs}
\title{Functional Equation for Theta Series  }

\author{Jae-Hyun Yang}
\address{Department of Mathematics, Inha University,
Incheon 402-751, Korea}
\email{jhyang@inha.ac.kr }

\begin{document}

\newtheorem{theorem}{Theorem}[section]
\newtheorem{lemma}{Lemma}[section]
\newtheorem{proposition}{Proposition}[section]
\newtheorem{remark}{Remark}[section]
\newtheorem{definition}{Definition}[section]

\renewcommand{\theequation}{\thesection.\arabic{equation}}
\renewcommand{\thetheorem}{\thesection.\arabic{theorem}}
\renewcommand{\thelemma}{\thesection.\arabic{lemma}}
\newcommand{\BR}{\mathbb R}
\newcommand{\BQ}{\mathbb Q}
\newcommand{\bn}{\bf n}
\def\charf {\mbox{{\text 1}\kern-.24em {\text l}}}
\newcommand{\BC}{\mathbb C}
\newcommand{\BZ}{\mathbb Z}



\begin{abstract}
{In this short paper, we find the transformation formula for the
theta series under the action of the Jacobi modular group on the
Siegel-Jacobi space. This formula generalizes the formula (5.1)
obtained by Mumford in \cite[p.\,189]{Mum1}.}
\end{abstract}
\maketitle

\newcommand\tr{\triangleright}
\newcommand\al{\alpha}
\newcommand\be{\beta}
\newcommand\g{\gamma}
\newcommand\gh{\Cal G^J}
\newcommand\G{\Gamma}
\newcommand\de{\delta}
\newcommand\e{\epsilon}
\newcommand\z{\zeta}
\newcommand\vth{\vartheta}
\newcommand\vp{\varphi}
\newcommand\om{\omega}
\newcommand\p{\pi}
\newcommand\la{\lambda}
\newcommand\lb{\lbrace}
\newcommand\lk{\lbrack}
\newcommand\rb{\rbrace}
\newcommand\rk{\rbrack}
\newcommand\s{\sigma}
\newcommand\w{\wedge}
\newcommand\fgj{{\frak g}^J}
\newcommand\lrt{\longrightarrow}
\newcommand\lmt{\longmapsto}
\newcommand\lmk{(\lambda,\mu,\kappa)}
\newcommand\Om{\Omega}
\newcommand\ka{\kappa}
\newcommand\ba{\backslash}
\newcommand\ph{\phi}
\newcommand\M{{\Cal M}}
\newcommand\bA{\bold A}
\newcommand\bH{\bold H}

\newcommand\Hom{\text{Hom}}
\newcommand\cP{\Cal P}
\newcommand\cH{\Cal H}

\newcommand\pa{\partial}

\newcommand\pis{\pi i \sigma}
\newcommand\sd{\,\,{\vartriangleright}\kern -1.0ex{<}\,}
\newcommand\wt{\widetilde}
\newcommand\fg{\frak g}
\newcommand\fk{\frak k}
\newcommand\fp{\frak p}
\newcommand\fs{\frak s}
\newcommand\fh{\frak h}
\newcommand\Cal{\mathcal}

\newcommand\fn{{\frak n}}
\newcommand\fa{{\frak a}}
\newcommand\fm{{\frak m}}
\newcommand\fq{{\frak q}}
\newcommand\CP{{\mathcal P}_g}
\newcommand\Hgh{{\mathbb H}_g \times {\mathbb C}^{(h,g)}}
\newcommand\BD{\mathbb D}
\newcommand\BH{\mathbb H}
\newcommand\CCF{{\mathcal F}_g}
\newcommand\CM{{\mathcal M}}
\newcommand\Ggh{\Gamma_{g,h}}
\newcommand\Chg{{\mathbb C}^{(h,g)}}
\newcommand\Yd{{{\partial}\over {\partial Y}}}
\newcommand\Vd{{{\partial}\over {\partial V}}}

\newcommand\Ys{Y^{\ast}}
\newcommand\Vs{V^{\ast}}
\newcommand\LO{L_{\Omega}}
\newcommand\fac{{\frak a}_{\mathbb C}^{\ast}}

%
%
\begin{section}{{\bf Introduction}}
\setcounter{equation}{0} For a given fixed positive integer $g$,
we let
$${\mathbb H}_g=\,\{\,\Omega\in \BC^{(g,g)}\,|\ \Om=\,^t\Om,\ \ \ \text{Im}\,\Om>0\,\}$$
be the Siegel upper half plane of degree $g$ and let
$$\G_g=\{ \g\in \BZ^{(2g,2g)}\ \vert \ ^t\!\g J_g\g= J_g\ \}$$
be the Siegel modular  group of degree $g$, where $F^{(k,l)}$
denotes the set of all $k\times l$ matrices with entries in a
commutative ring $F$ for two positive integers $k$ and $l$,
$^t\!M$ denotes the transpose matrix of a matrix $M,\
\text{Im}\,\Om$ denotes the imaginary part of $\Om$ and
$$J_g=\begin{pmatrix} 0&I_g\\
                   -I_g&0\end{pmatrix}.$$

For two positive integers $g$ and $m$, we consider the Heisenberg
group
$$H_{\BZ}^{(g,m)}:=\big\{ \,(\lambda,\mu;\kappa)\,|\
\lambda,\,\mu\in\BZ^{(m,g)},\ \kappa\in\BZ^{(m,m)},\
\kappa+\mu\,{}^t\!\lambda\ \,symmetric\,\big\}$$

\noindent endowed with the following multiplication law
$$(\lambda,\mu;\kappa)\circ
(\lambda',\mu';\kappa'):=(\lambda+\lambda',\mu+\mu';\kappa+\kappa'+\lambda\,{}^t\!\mu'-\mu\,{}^t\!\lambda').$$

\noindent We let
\begin{equation*}
\G_{g,m}:=\G_g\ltimes H_\BZ^{(g,m)}\quad ( \textrm{semi-direct
proudct)}
\end{equation*}
be the Jacobi modular group endowed with the following
multiplication law
$$
(\g,(\lambda,\mu;\kappa))\cdot(\g',(\lambda',\mu';\kappa')) =\,
\big(\g\g',(\widetilde{\lambda}+\lambda',\widetilde{\mu}+
\mu';\kappa+\kappa'+{\widetilde
\lambda}\,{}^t\!\mu'-{\widetilde\mu}\,{}^t\!\lambda' )\big)$$ with
$\g,\g'\in \G_g,\ \lambda,\lambda',\mu,\mu'\in\BZ^{(m,g)},\
\kappa,\,\kappa'\in\BZ^{(m,m)}$ and
$(\widetilde{\lambda},\widetilde{\mu})=(\lambda,\mu)\g'$. Then
$\G_{g,m}$ acts on the Siegel-Jacobi space $\BH_{g,m}:=\BH_g\times
\BC^{(m,g)}$ properly discontinuously by
\begin{equation}
(\g,(\lambda,\mu;\kappa))\cdot (\Om,Z)=\big(
(A\Om+B)(C\Om+D)^{-1},(Z+\lambda \Om+\mu)
(C\Omega+D)^{-1}), \end{equation} where $\g=\begin{pmatrix} A&B\\
C&D\end{pmatrix} \in \G_g,\ \lambda,\mu\in\BZ^{(m,g)},\
\kappa\in\BZ^{(m,m)}$ and $(\Om,Z)\in
\BH_{g,m}$\,(cf.\,\cite{YJ1},\,\cite{YJ2},\,\cite{YJ4},\,\cite{YJ5}).
A fundamental domain for $\G_{g,m}\backslash \BH_{g,m}$ was found
by the author in \cite{YJ3}. Let $\G_{\vartheta,g}$ be the theta
group consisting of all element
$\g=\begin{pmatrix} A&B\\
C&D\end{pmatrix} \in \G_g$ such that the diagonal entries of
matrices ${}^t\! AC$ and ${}^tBD$ are {\it even} integers. We set
$$ \G_{\vartheta,g,m}:=\G_{\vartheta,g}\ltimes
H_\BZ^{(g,m)}.$$ We consider the theta series
\begin{equation}
\Theta(\Om,Z):=\sum_{A\in\BZ^{(m,g)}} e^{\pi i
\,\s(A\Om\,{}^t\!A+2A\,{}^t\!Z)},\quad (\Omega,Z)\in\BH_{g,m}.
\end{equation}
Here $\s(T)$ denotes the trace of a square matrix $T$.

\vspace{0.2in} In \cite[p.\,189]{Mum1}, Mumford considered the
case $m=1$ and proved the following functional equation
\begin{eqnarray}
&&\Theta\big((A\Om+B)(C\Om+D)^{-1},Z(C\Om+D)^{-1}\big)\\
&=& \zeta(\g)\,e^{\pi i\{\,Z(C\Om+D)^{-1}C\,{}^tZ\,\}}\,
\det(C\Om+D)^{1/2}\,\Theta(\Om,Z),\notag
\end{eqnarray}
where $\g=\begin{pmatrix} A&B\\
C&D\end{pmatrix} \in \G_{\vartheta,g}$ and $\zeta(\g)$ is an
eighth root of $1$. \vspace{0.1in} In this short article, we
consider the case of an arbitrary positive integer $m$ and then
prove the following functional equation.

\begin{theorem}
For any
${\tilde\g}=(\g,(\lambda,\mu;\kappa))\in\G_{\vartheta,g,m}$ with
$\g=\begin{pmatrix} A&B\\
C&D\end{pmatrix} \in \G_{\vartheta,g}$, we obtain the following
functional equation
\begin{eqnarray}
&&\Theta\big((A\Om+B)(C\Om+D)^{-1},(Z+\lambda \Om+\mu)(C\Om+D)^{-1}\big)\\
&=& \zeta(\tilde\g)\,e^{\pi
i\,\s\{\,(Z+\lambda\Om+\mu)(C\Om+D)^{-1}C\,{}^t(Z+\lambda\Om+\mu)-\lambda
\Om\,{}^t\!\lambda-2\lambda\,{}^t\,Z\} }\, \det(C\Om+D)^{\frac
m2}\,\Theta(\Om,Z),\notag
\end{eqnarray}
where $\zeta({\tilde\g})$ is an eighth root of $1$.
\end{theorem}
We observe that the formula (1.4) generalizes the formula (1.3)
with $m=1$ and $\lambda=\mu=0$. For a positive integer $N$, we put
$$\G_0(N)=\left\{ \begin{pmatrix} a&b\\
                   c&d\end{pmatrix}\in SL(2,\BZ)\,\big|\ c\equiv 0\,(
                   \textrm{mod}\,N)\,\right\}$$
and
$$\theta(\tau)=\sum_{r=-\infty}^\infty e^{2\pi i\,r^2\tau},\quad
\tau \in\BH_1.$$

 In \cite[(Werke)\,pp.\,939--940]{H}, Hecke showed that

\begin{equation}
\theta\big( (a\tau+b)(c\tau+d)^{-1}\big)=\epsilon_d^{-1}\left(
{c\over d}\right)\,(c\tau+d)^{1/2}\,\theta(\tau),\quad
\begin{pmatrix} a&b\\
                   c&d\end{pmatrix}\in \G_0(4),
\end{equation}
where $\epsilon_d=1$ or $i$ according to $d\equiv 1$ or $3$\,($
\textrm{mod}\,4)$ and $\left( {c\over d}\right)$ denotes the
quadratic residue symbol\,(cf.\,\cite[p.\,442]{S}).

\vskip 0.2cm \noindent {\bf Notations\,:} \ \ We denote by $\BZ$
and $\BC$ the ring of integers, and the field of complex numbers
respectively. $\BC^{\times}$ denotes the multiplicative group of
nonzero complex numbers. The symbol ``:='' means that the
expression on the right is the definition of that on the left. For
two positive integers $k$ and $l$, $F^{(k,l)}$ denotes the set of
all $k\times l$ matrices with entries in a commutative ring $F$.
For a square matrix $A\in F^{(k,k)}$ of degree $k$, $\sigma(A)$
denotes the trace of $A$. For any $M\in F^{(k,l)},\ ^t\!M$ denotes
the transpose matrix of $M$. $I_n$ denotes the identity matrix of
degree $n$. We put $i=\sqrt{-1}.$ For $z\in\BC,$ we define
$z^{1/2}=\sqrt{z}$ so that $-\pi / 2 < \ \arg (z^{1/2})\leqq
\pi/2.$ Further we put $z^{\kappa/2}=\big(z^{1/2}\big)^\kappa$ for
every $\kappa\in\BZ.$

\end{section}
\vskip 0.3cm
%
%
\begin{section}{{\bf Proof of Theorem 1.1}}
\setcounter{equation}{0}

Let $\tilde\g=(\g,(\la,\mu;\kappa))$ be an element of $\G_{g,m}$
with $\g=\begin{pmatrix} A&B\\ C&D
\end{pmatrix}\in \G_g$ and $(\Om,Z)\in\BH_{g,m}$ with $\Om\in\BH_g$ and $Z\in\BC^{(m,g)}.$
If we put $(\Om_*,Z_*):=\tilde\g\cdot(\Om,Z),$ then we have
\begin{eqnarray*}
&\Om_*=\g\cdot \Om=(A\Om+B)(C\Om+D)^{-1},\\
&Z_*=(Z+\la \Om+\mu)(C\Om+D)^{-1}.
\end{eqnarray*}

First of all we shall show that if the formula (1.4) holds for
$\tilde\g_1,\,\tilde\g_2\in \G_{g,m},$ then it hold for
$\tilde\g_1\tilde\g_2.$ To prove this fact, we consider the
function $J:\G_{g,m}\times \BH_{g,m}\lrt\BC^{\times}$ defined by
\begin{equation*}
J\big(\tilde\g, (\Om,Z)\big):=\,e^{\pi
i\,\s\{\,(Z+\lambda\Om+\mu)(C\Om+D)^{-1}C\,{}^t(Z+\lambda\Om+\mu)-\lambda
\Om\,{}^t\!\lambda-2\lambda\,{}^t\!Z-\kappa-\mu\,{}^t\!\lambda\}
},
\end{equation*}
where
$\tilde\g=(\g,(\la,\mu;\kappa))\in\G_{g,m}$ with $\g=\begin{pmatrix} A&B\\
C&D\end{pmatrix} \in \G_g,\ \lambda,\mu\in\BZ^{(m,g)},\
\kappa\in\BZ^{(m,m)}$ and $(\Om,Z)\in \BH_{g,m}$. By a direct
computation or a geometrical method\,(cf,\,\cite[p.\,1332]{YJ2}),
we can show that $J$ is an automorphic factor for $\G_{g,m}$ on
$\BH_{g,m}$, that is, it satisfies the following relation
\begin{equation*}
J\big(\tilde\g_1\tilde\g_2, (\Om,Z)\big)=
J\big(\tilde\g_1,\tilde\g_2\!\cdot\!
(\Om,Z)\big)\,J\big(\tilde\g_2, (\Om,Z)\big)
\end{equation*}
for any $\tilde\g_1,\,\tilde\g_2\in \G_{g,m}$ and
$(\Om,Z)\in\BH_{g,m}.$ It is easy to see that the map $J_*:
\G_{g,m}\times \BH_{g,m}\lrt\BC^{\times}$ defined by
\begin{equation*}
J_*\big(\tilde\g, (\Om,Z)\big):=J\big(\tilde\g,
(\Om,Z)\big)\!\cdot\!\det(C\Om+D)^{\frac m2}
\end{equation*}
is an automorphic factor for $\G_{g,m}$ on $\BH_{g,m}$, where
$\tilde\g=(\g,(\la,\mu;\kappa))\in\G_{g,m}$ with $\g=\begin{pmatrix} A&B\\
C&D\end{pmatrix} \in \G_g$ and $(\Om,Z)\in \BH_{g,m}$. It is
easily seen that $J_*\big(\tilde\g, (\Om,Z)\big)$ can be written
as
\begin{equation*}
J_*\big(\tilde\g, (\Om,Z)\big)=e^{-\pi
i\,\s(\kappa+\mu\,{}^t\!\lambda)}\cdot e^{\pi
i\,\s\{\,(Z+\lambda\Om+\mu)(C\Om+D)^{-1}C\,{}^t(Z+\lambda\Om+\mu)-\lambda
\Om\,{}^t\!\lambda-2\lambda\,{}^t\!Z\} }\,\det(C\Om+D)^{\frac m2}.
\end{equation*}

\noindent We observe that $e^{-\pi
i\,\s(\kappa+\mu\,{}^t\!\lambda)}=\pm 1$ because
$\s(\kappa+\mu\,{}^t\!\lambda)$ is an integer. Thus we see that if
the formula (1.4) holds for $\tilde\g_1,\,\tilde\g_2\in \G_{g,m},$
then it hold for $\tilde\g_1\tilde\g_2.$

\vskip0.2cm We recall
(cf.\,\cite[p.\,326]{F},\,\cite[p.\,210]{Mum1}) that $\G_g$ is
generated by the following elements
\begin{eqnarray*}
t_0(B):&=&\begin{pmatrix} I_g& B\\
                   0& I_g\end{pmatrix}\ \textrm{with any}\
                   B=\,{}^tB\in \BZ^{(g,g)},\\
g_0(\alpha):&=&\begin{pmatrix} {}^t\alpha & 0\\
                   0& \alpha^{-1}  \end{pmatrix}\ \textrm{with
                   any}\ \alpha\in GL(g,\BZ),\\
-J_g:&=&\begin{pmatrix} 0& -I_g\\
                   I_g&\ 0\end{pmatrix}.
\end{eqnarray*}

\noindent Obviously the following matrices
\begin{eqnarray*}
t_e(B):&=&\begin{pmatrix} I_g& B\\
                   0& I_g\end{pmatrix}\ \textrm{with any}\
                   B=\,{}^tB\in \BZ^{(g,g)}\ even\ \textrm{diagonals},\\
g_0(\alpha):&=&\begin{pmatrix} {}^t\alpha & 0\\
                   0& \alpha^{-1}  \end{pmatrix}\ \textrm{with
                   any}\ \alpha\in GL(g,\BZ),\\
-J_g:&=&\begin{pmatrix} 0& -I_g\\
                   I_g&\ 0\end{pmatrix}.
\end{eqnarray*}
generate the theta group $\G_{\vartheta,g}$. Therefore the
following elements $s(\lambda,\mu;\kappa),\ t(B),\,g(\alpha)$ and
$\s_g$ of $\G_{\vartheta,g,m}$ defined by
\begin{eqnarray*}
&& s(\la,\mu;\kappa)=\big( I_{2g},(\la,\mu;\kappa)\big)\ \textrm{with}\ \la,\mu\in
\BZ^{(m,g)}\ \textrm{and}\ \kappa\in\BZ^{(m,m)} ,\\
&&t(B)=\left(\begin{pmatrix}
I_g&B\\0&I_g\end{pmatrix},(0,0;0)\right)\ \textrm{with any}\
                   B=\,{}^tB\in \BZ^{(g,g)}\ even\ \textrm{diagonals},\\
&& g(\alpha)=\left( \begin{pmatrix} {}^t\alpha & 0\\
                   0& \alpha^{-1}  \end{pmatrix},(0,0;0)\right)\
\textrm{with}\ \alpha\in GL(g,\BZ),\\
 &&\s_g=\left(\begin{pmatrix}
0&-I_g\\I_g&0\end{pmatrix},(0,0;0)\right)
\end{eqnarray*}
generate the group $\G_{\vartheta,g,m}.$

\par
\vskip3mm \noindent {\bf Case I.} $\tilde\g=s(\la,\mu;\kappa)$
with $\la,\mu\in\BZ^{(m,g)} \ \textrm{and}\ \kappa\in\BZ^{(m,m)}.$
\vskip 0.2cm In this case, we have
$$\Om_*=\Om\quad \textrm{and}\quad Z_*=Z+\la \Om+\mu.$$
Then we have
\begin{eqnarray*}
& & \Theta(\Om,Z+\la \Om+\mu)\\
&=&\sum_{A\in\BZ^{(m,g)} } e^{\pi i\,\s\{
\,A\Om\,{}^t\!A+\,2\,A\,{}^t\!(Z+\la
\Om+\mu)\,\} }
\\
&=&\,e^{-\pi i
\,\s(\la\Om\,{}^t\!\la+2\,\la\,{}^t\!Z)}\,\sum_{A\in\BZ^{(m,g)}}
e^{\pi i\,\s\big\{\,(A+\la)\,\Om\,{}^t\!(A+\la)
+\,2\,(A+\la)\,{}^t\!Z\,\big\} }
\\
&=&\,e^{-\pi i
\,\s(\la\Om\,{}^t\!\la+2\,\la\,{}^t\!Z)}\,\Theta(\Om,Z).
\end{eqnarray*}
Here we may take $\zeta(\tilde\g)=1.$ Therefore this proves the
formula (1.4) in the case $\tilde\g=s(\la,\mu;\kappa)$.

\par
\vskip3mm\noindent {\bf Case II.} $\tilde\g=t(B)$ with
$B=\,{}^tB\in\BZ^{(g,g)}$ even diagonal. \vskip 0.2cm In this
case, we have
$$\Om_*=\Om+B\quad \textrm{and}\quad Z_*=Z.$$
Then we have
\begin{eqnarray*}
&& \Theta(\Om+B,Z)\\
&=&\sum_{A\in\BZ^{(m,g)}} e^{\pi i\,\s \{ A(\Om+B)\,{}^t\!A+2A\,{}^t\!Z\,\} }\\
&=& \sum_{A\in\BZ^{(m,g)}} e^{\pi i\,\s ( A\Om\,{}^t\!A+2A\,{}^t\!Z\,) }\cdot
e^{\pi i\,\s(AB\,{}^t\!A)}    \\
&=&\sum_{A\in\BZ^{(m,g)}} e^{\pi i\,\s (
A\Om\,{}^t\!A+2A\,{}^t\!Z\,) }\quad \quad ( \textrm{because}\
\s(AB\,{}^t\!A)\in 2\,\BZ)\\
&=&\Theta(\Om,Z)
\end{eqnarray*}
Here we note that $\s(AB\,{}^t\!A)\in 2\,\BZ$ because the diagonal
entries of $B$ is even integers. Now we may take
$\zeta(\tilde\g)=1.$ Therefore this proves the formula (1.4) in
the case $\tilde\g=t(B)$.

\par
\vskip3mm \noindent {\bf Case III.} $\tilde\g=g(\alpha)=\left( \begin{pmatrix} {}^t\alpha & 0\\
                   0& \alpha^{-1}  \end{pmatrix},(0,0;0)\right)\
\textrm{with}\ \alpha\in GL(g,\BZ).$ \vskip 0.2cm In this case, we
have
$$\Om_*=\,{}^t\alpha\Om\,\alpha\quad \textrm{and}\quad Z_*=Z\alpha.$$
Then we obtain
\begin{eqnarray*}
&& \Theta(\,{}^t\alpha\Om\,\alpha,Z\alpha )\\
&=&\sum_{A\in\BZ^{(m,g)}} e^{\pi i\,\s \{ A(\,{}^t\alpha\Om\,\alpha)\,{}^t\!A\,+\,2\,A\,{}^t(Z\,\alpha)\} }\\
&=& \sum_{A\in\BZ^{(m,g)}} e^{\pi i\,\s
\{\,(A\,{}^t\alpha)\,\Om\,{}^t(A\,{}^t\alpha)\,+\, 2\,(A\,{}^t\alpha)Z\,\} }    \\
&=&\Theta(\Om,Z).
\end{eqnarray*}
We observe that the formula (1.4) reduces to the formula
\begin{equation}
\Theta\big(\,{}^t\alpha\Om\,\alpha,Z\alpha\big)=\,\zeta(\tilde\g)\,\big(\det
\alpha^{-1}\big)^{m/2}\,\Theta(\Om,Z).
\end{equation}
If we take $\zeta(\tilde\g)=\big(\det \alpha\big)^{m/2},$ the
formula (2.1) coincides with $\Theta(\Om,Z).$ Since $\det
\alpha=\pm 1,\ \zeta(\tilde\g)$ is a fourth root of $1$.
Therefore this proves the formula (1.4) in the case
$\tilde\g=g(\alpha)$ with $\alpha\in GL(g,\BZ)$.

\par
\vskip3mm \noindent {\bf Case IV.}
$\tilde\g=\s_g=\left(\begin{pmatrix}
0&-I_g\\I_g&0\end{pmatrix},(0,0;0)\right).$ \vskip 0.2cm In this
case, we have
$$\Om_*=-\Om^{-1}\quad \textrm{and}\quad Z_*=Z\Om^{-1}.$$
We can prove the formula (1.4) using the Poisson Summation
Formula. \vskip 0.3cm\noindent {\bf Lemma 2.1.} For a fixed
element $(\Om,Z)\in {\mathbb H}_{g,m},$ we obtain the following
\begin{equation*}
\int_{\BR^{(m,g)}} e^{\pi i\, \sigma ( x\Om\, {}^tx + 2x\,{}^tZ )}
dx_{11}\cdots dx_{mg}=\left( \det \left( { {\Omega}\over
i}\right)\right)^{-{\frac m2}} e^{-\pi i\,
\sigma(Z\Om^{-1}\,{}^tZ)},
\end{equation*}
where $x=(x_{ij})\in\BR^{(m,g)}.$ \vskip 0.3cm\noindent {\it
Proof.} By a simple computation, we see that
$$e^{\pi i\, \sigma ( x\Om\, {}^tx +
2x\,{}^tZ )}= e^{-\pi i\,\sigma (Z\Om^{-1}\,{}^tZ )}\cdot e^{\pi
i\,\sigma \{(x+Z\Om^{-1})\Om\,{}^t(x+Z\Om^{-1})\} }.$$ Since the
real Jacobi group $Sp(g,\BR)\ltimes H_\BR^{(m,g)}$ acts on
${\mathbb H}_{g,m}$ holomorphically, we may put
$$\Om=\,i\,A\,{}^t\!A,\quad Z=iV,\quad\  A\in\BR^{(g,g)},\quad
V=(v_{ij})\in\BR^{(m,g)}.$$

\begin{eqnarray*}
& & \int_{\BR^{(m,g)}}  e^{\pi i\, \sigma ( x\Om\, {}^tx +
2x\,{}^tZ )} dx_{11}\cdots dx_{mg} \\
&=& e^{-\pi i\,\sigma (Z\Omega^{-1}\,{}^tZ)} \int_{\BR^{(m,g)}}
e^{\pi i\,\sigma [\{
x+iV(iA\,{}^t\!A)^{-1}\}(iA\,{}^t\!A)\,{}^t\!\{
x+iV(iA\,{}^t\!A)^{-1}\} ]}\,dx_{11}\cdots dx_{mg}\\
&=&e^{-\pi i\,\sigma (Z\Omega^{-1}\,{}^tZ)} \int_{\BR^{(m,g)}}
e^{\pi i\,\sigma [\{ x+V(A\,{}^t\!A)^{-1}\}A\,{}^t\!A\,{}^t\!\{
x+V(A\,{}^t\!A)^{-1}\} ]}\,dx_{11}\cdots dx_{mg}\\
&=& e^{-\pi i\,\sigma (Z\Omega^{-1}\,{}^tZ)} \int_{\BR^{(m,g)}}
e^{-\pi \,\sigma\{ (uA)\,{}^t\!(uA)\} }\,du_{11}\cdots
du_{mg}\quad \big(\,{\rm Put}\ u=
x+V(A\,{}^t\!A)^{-1}=(u_{ij}) \,\big)\\
&=& e^{-\pi i\,\sigma (Z\Omega^{-1}\,{}^tZ)} \int_{\BR^{(m,g)}}
e^{-\pi \,\sigma (w\,{}^t\!w)} (\det A)^{-m}\,dw_{11}\cdots
dw_{mg}\quad \big(\,{\rm Put}\ w=uA=(w_{ij})\,\big)\\
&=& e^{-\pi i\,\sigma (Z\Omega^{-1}\,{}^tZ)} \, (\det A)^{-m}\cdot
\left( \prod_{i=1}^m \prod_{j=1}^g \int_\BR e^{-\pi\,
w_{ij}^2}\,dw_{ij}\right)\\
&=& e^{-\pi i\,\sigma (Z\Omega^{-1}\,{}^tZ)} \, (\det A)^{-m}\quad
\big(\,{\rm because}\ \int_\BR e^{-\pi\,
w_{ij}^2}\,dw_{ij}=1\quad {\rm for\ all}\ i,j\,\big)\\
&=& e^{-\pi i\,\sigma (Z\Omega^{-1}\,{}^tZ)} \, \left( \det \big(
A\, {}^t\!A \big)\right)^{-{\frac m2}}\\
&=& e^{-\pi i\,\sigma (Z\Omega^{-1}\,{}^tZ)} \, \left( \det \left(
{ {\Omega}\over i } \right)\right)^{-{\frac m2}}.
\end{eqnarray*}

\noindent This completes the proof of Lemma 2.1. \hfill $\square$

\vskip 0.2cm For an element $(\Om,Z)\in {\mathbb H}_{g,m},$ we
define the function $f_{\Omega,Z}$ on $\BR^{(m,g)}$ by
\begin{equation}
f_{\Omega,Z}(x):=e^{\pi i\, \sigma ( x\Om\, {}^tx + 2x\,{}^tZ
)},\quad x=(x_{ij})\in\BR^{(m,g)}.
\end{equation}

\noindent By the Poisson summation formula, we obtain
$$ \sum_{A\in\BZ^{(m,g)}} f_{\Omega,Z}(A) =  \sum_{A\in\BZ^{(m,g)}} {\widehat
f}_{\Omega,Z}(A),$$ where ${\widehat f}_{\Omega,Z}$ is the Fourier
transform of $f_{\Om,Z}$ given by
$${\widehat
f}_{\Omega,Z}(y)=\int_{\BR^{(m,g)}} f_{\Om,Z}(x)\,e^{2\pi
i\,\sigma (\,{}^t\!x y)}\,dx_{11}\cdots dx_{mg}.$$ Then we have

\begin{eqnarray*}
\Theta(\Omega,Z)&=&\sum_{A\in\BZ^{(m,g)}} {\widehat
f}_{\Omega,Z}(A)\\
&=&\sum_{A\in\BZ^{(m,g)}}\int_{\BR^{(m,g)}}  f_{\Om,Z}(x)\,e^{2\pi
i\,\sigma
(\,{}^t\!xA)}\,dx_{11}\cdots dx_{mg}\\
&=&\sum_{A\in\BZ^{(m,g)}}\int_{\BR^{(m,g)}} e^{\pi i\, \sigma (
x\Om\, {}^tx + 2x\,{}^tZ )} \,e^{2\pi i\,\sigma
(\,{}^t\!xA)}\,dx_{11}\cdots dx_{mg}\\
&=&\sum_{A\in\BZ^{(m,g)}}\int_{\BR^{(m,g)}} e^{\pi i\, \sigma \{
x\Om\, {}^tx + 2x\,{}^t(Z+A)\} } \,dx_{11}\cdots dx_{mg}\\
&=&\sum_{A\in\BZ^{(m,g)}} \left( \det \left( { {\Omega}\over i }
\right)\right)^{-{\frac m2}}\,e^{-\pi i\,\sigma \{
(Z+A)\Omega^{-1}\,{}^t(Z+A)\} }\quad (\,{\rm by\ Lemma\ 2.1\,})\\
&=& \left( \det \left( { {\Omega}\over i } \right)\right)^{-{\frac
m2}}\,\sum_{A\in \BZ^{(m,g)}} e^{-\pi i\,\sigma (
Z\Omega^{-1}\,{}^t\!Z+
A\Omega^{-1}\,{}^t\!A+2A\Omega^{-1}\,{}^t\!Z)}\\
&=& \left( \det \left( { {\Omega}\over i } \right)\right)^{-{\frac
m2}} \, e^{-\pi i\,\sigma ( Z\Omega^{-1}\,{}^t\!Z)}   \,\sum_{A\in
\BZ^{(m,g)}} e^{\pi i \,\sigma
\{(-A)(-\Om^{-1})\,{}^t\!(-A)+2(-A)\,{}^t\!(Z\Om^{-1})\} }\\
&=& \left( \det \left( { {\Omega}\over i } \right)\right)^{-{\frac
m2}} \, e^{-\pi i\,\sigma ( Z\Omega^{-1}\,{}^t\!Z)}
\,\Theta(-\Om^{-1},Z\Om^{-1}).
\end{eqnarray*}
Therefore we obtain the formula
\begin{equation}
\Theta\big(-\Om^{-1},Z\Omega^{-1}\big)=\,e^{\pi i\,\s(
Z\Om^{-1}\,{}^t\!Z)}\,\left( \det \left( {\Om\over
i}\right)\right)^{m/2}\,\Theta(\Om,Z).
\end{equation}
The fact that $\zeta(\tilde\g)^8=1$ follows from the formula
(2.3). Indeed we may take $\zeta(\tilde\g)=\det\left( {{I_g}\over
i} \right)^{m/2}.$ Therefore this proves the formula (1.4) in the
case $\tilde\g=\s_g$. Finally we complete the proof of Theorem
1.1.

\vspace{0.1in}\noindent {\bf Remark 2.1.} Let $m$ be an {\it odd}
positive integer. According to the formula (1.4), we see that
$\Theta(\Om,0)$ is a modular form of half integral weight ${\frac
m2}$ with respect to
$\G_{\vartheta,g}$\,(cf.\,\cite[p.\,200]{Mum1},\,\cite{S}). We may
say that the theta series $\Theta(\Om,Z)$ is a Jacobi form of half
integral weight ${\frac m2}$ and index $I_m$ with respect to
$\G_{\vartheta,g}$\,(cf.\,\cite{YJ1},\,\cite{YJ2}). This means
that $\Theta(\Om,Z)$ may be regarded as an automorphic form on a
two-fold covering of the Jacobi group \,(cf.\,\cite{Sat}). Indeed
the theta series $\Theta(\Om,Z)$ is closely related to the Weil
representation of the Jacobi group\,(cf.\,\cite{W},\,\cite{YJH}).
The function $f_{\Omega,Z}$ is a covariant map for the
Weil-Schr{\"o}dinger representation\,(cf.\,\cite{YJH}).

\vspace{0.1in}\noindent {\bf Remark 2.2.} Olav K. Richter \cite{R}
obtained the transformation formula for theta functions that is
more general than the formula (1.4). It is my pleasure to thank
him for letting me know his paper \cite{R}. But our proof is quite
different from his. In fact, our formula (1.4) is a combination of
the transformation laws (2) and (3) in \cite{R}.

\end{section}


\vskip 0.1cm








\vspace{0.5cm}



\begin{thebibliography}{99}


\bibitem{F} E. Freitag, {\em Siegelsche Modulfunktionen}, Grundlehren de mathematischen Wissenschaften {\bf 55},
Springer-Verlag, Berlin-Heidelberg-New York (1983).

\bibitem{H} E. Hecke, {\em Herleitung des Euler-Produktes der Zetafunktion und einiger L-Reihnen
aus ihrer Funktionalgleichung}, Math. Ann. {\bf 119}\,(1944),
266-287 (=Werke, 919-940).


\bibitem{Mum1} D. Mumford, {\em Tata Lectures on Theta I,} Progress
in Math. {\bf 28}, Boston-Basel-Stuttgart (1983).


\bibitem{R} O. K. Richter, {\em On Transformation Laws for Theta Functions,}
Rocky Mountain J. of Math., Vol. {\bf 34}, No. 4 (2004),
1473-1481.



\bibitem{Sat} I. Satake, {\em Fock representations and theta functions},
Ann. Math. Study {\bf 66}\,(1969), 393--405.


\bibitem{S} G. Shimura, {\em On modular forms of half integral
weight }, Ann. of Math., {\bf 97}\,(1973), 440-481; Collected
Papers, 1967-1977, Vol.\,II, Springer-Verlag\,(2002), 532-573.

\bibitem{W} A. Weil, {\em Sur certains groupes d'operateurs unitares},
Acta Math., {\bf 111}\,(1964), 143--211; Collected
Papers\,(1964-1978), Vol.\,III, Springer-Verlag\,(1979), 1-69.


\bibitem{YJ1} J.-H. Yang, \textit{Singular Jacobi forms,}
Trans. of American Math. Soc. {\bf 347}, No. 6 (1995), 2041-2049.

\bibitem{YJ2} J.-H. Yang, \textit{Construction of vector valued
modular forms from Jacobi forms,} Canadian J. of Math. {\bf 47\,
(6)} (1995), 1329-1339.

\bibitem{YJ3} J.-H. Yang, \textit{A note on a fundamental domain for Siegel-Jacobi space,} Houston Journal of
Mathematics, Vol. $ \textbf{32}$, No. 3 (2006), 701--712.

\bibitem{YJ4} J.-H. Yang, {\em Inavraint metrics and Laplacians on Siegel-Jacobi space,}
Journal of Number Theory, Vol. $ \textbf{127}$ (2007), 83-102.

\bibitem{YJ5} J.-H. Yang, {\em A partial Cayley transform of Siegel-Jacobi
disk,} J. Korean Math. Soc. {\bf 45}, No. 3 (2008), 781-794.

\bibitem{YJH} J.-H. Yang, {\em Theta series  associated with the Weil-Schr{\"o}dinger representation,}
arXiv:0709.007v1 [math.NT] (2007).












\end{thebibliography}
\end{document}